\documentclass[11pt]{amsart}
\usepackage{amsmath}
\usepackage{amssymb}
\usepackage{amsthm}
\usepackage{amscd}
\usepackage{enumerate}
\usepackage{enumitem}
\usepackage{verbatim}
\usepackage[all]{xy}
\usepackage{mathrsfs}
\usepackage{mathabx}
\usepackage{url}
\usepackage{hyperref}

\theoremstyle{plain}
\newtheorem{thm}{Theorem}[section]
\newtheorem{prop}[thm]{Proposition}
\newtheorem{lem}[thm]{Lemma}
\newtheorem{cor}[thm]{Corollary}

\theoremstyle{definition}
\newtheorem{dfn}[thm]{Definition}

\newcommand{\diag}{\mathrm{diag}}

\newcommand{\Fr}{\mathrm{Fr}}

\newcommand{\Hom}{\mathrm{Hom}}

\newcommand{\id}{\mathrm{id}}

\newcommand{\rank}{\mathrm{rank}}

\newcommand{\Supp}{\mathrm{Supp}}
\newcommand{\Sym}{\mathrm{Sym}}

\newcommand{\cB}{\mathcal{B}}

\newcommand{\cE}{\mathcal{E}}

\newcommand{\cG}{\mathcal{G}}
\newcommand{\cH}{\mathcal{H}}

\newcommand{\cP}{\mathcal{P}}

\newcommand{\cU}{\mathcal{U}}
\newcommand{\cV}{\mathcal{V}}

\newcommand{\vep}{\varepsilon}

\newcommand{\bF}{\mathbb{F}}

\newcommand{\bZ}{\mathbb{Z}}

\newcommand{\frn}{\mathfrak{n}}

\newcommand{\Cinf}{\mathbb{C}_{\infty}}

\makeatletter
    
    \@addtoreset{equation}{section}
\makeatother

\makeatletter
\renewcommand{\p@enumii}{}
\makeatother

\begin{document}

\title[Half-integrality and an upper bound for slopes]{Half-integrality and an upper bound for slopes of Drinfeld cuspforms of level $
\Gamma_1(t^r)$}

\author{Shin Hattori}
\address{Department of Natural Sciences, Tokyo City University, 1-28-1 Tamazutsumi, Setagaya-ku, Tokyo 158-8557, Japan}

\date{\today}


\begin{abstract}
Let $p$ be a rational prime, let $q>1$ be a power of $p$ and let $\bF_q$ be the field of $q$ elements. 
Let $k\geq 2$ and $r\geq 1$ be integers. Let $S_k(\Gamma_1(t^r))$ be the $\Cinf$-vector space of Drinfeld cuspforms of rank two, 
level $\Gamma_1(t^r)$ and weight $k$ over $
\bF_q(t)$.
In this paper, we prove that every finite $t$-adic slope of the Hecke operator $U_t$ on $S_k(\Gamma_1(t^r))$ belongs to $\frac{1}{2}\bZ$ and is at most $k-2$. 
If 
$q$ is odd, then every such slope is an 
integer.
\end{abstract}

\maketitle



\section{Introduction}\label{SecIntro}

Let $p$ be a rational prime and let $q>1$ be a power of $p$. Let $\bF_q$ be the field of $q$ elements and let $A=\bF_q[t]$ be the polynomial ring over $
\bF_q$.
Let $F=\bF_q(t)$, $F_\infty=\bF_q((1/t))$ and let $\Cinf$ be the $(1/t)$-adic completion of an algebraic closure of $F_\infty$. Let $\Omega=\Cinf\setminus 
F_\infty$ be the 
Drinfeld upper half plane. Let $r\geq 1$ be an integer and let
\[
\Gamma_1(t^r)=\left\{g\in\mathit{SL}_2(A)\ \middle|\ g\equiv \begin{pmatrix}
	1&*\\0&1
\end{pmatrix}\bmod t^r\right\}.
\]
Let $k\geq 2$ be an integer. We denote by $S_k(\Gamma_1(t^r))$ the $\Cinf$-vector space of Drinfeld cuspforms of rank two, level $\Gamma_1(t^r)$ and weight $k$.

Define the $U_t$ operator on $S_k(\Gamma_1(t^r))$ by
\[
(U_tf)(z)=\frac{1}{t}\sum_{b\in \bF_q}f\left(\frac{z+b}{t}\right)
\]
for any $f\in S_k(\Gamma_1(t^r))$ and $z\in \Omega$. 
It is known that every eigenvalue of $U_t$ is integral over $A$. We refer to the normalized $t$-adic additive valuation of an eigenvalue of $U_t$ as 
a ($t$-adic) slope of $U_t$.
Slopes are finite for nonzero eigenvalues, and the slope for the eigenvalue zero is infinite.

In this paper, we prove the following folklore conjectures on $t$-adic slopes in $S_k(\Gamma_1(t^r))$ (see, for example, \cite[\S4]{Ha_GVt}).

\begin{thm}[Theorem \ref{ThmMain}]\label{ThmMainIntro}
	Let $\alpha$ be any finite slope of $U_t$ acting on $S_k(\Gamma_1(t^r))$. Then the following hold.
	\begin{enumerate}
		\item\label{ThmMainIntro_HI} $\alpha\in \frac{1}{2}\bZ$. If $q$ is odd, then $\alpha\in \bZ$.
		\item\label{ThmMainIntro_Bound} $\alpha\leq k-2$.
	\end{enumerate}
\end{thm}

Note that, for any $p$-adic slope $\beta$ of elliptic cuspforms of level $\Gamma_1(p)$ and weight $k\geq 2$, 
the inequality $\beta\leq k-1$ follows from
\cite[Theorem 4.6.17]{Miy} for newforms, and from the relation $U^2-T_pU+p^{k-1}\langle p\rangle=0$ satisfied by the $U_p$ operator on the space of 
oldforms. However, neither proof extends directly to the Drinfeld case.
In particular, for oldforms, the corresponding relation for the Drinfeld case is $U^2-T_t U=0$ and this relation does not by itself yield the desired upper bound. 
This difficulty may be related to the lack of a suitable notion of weight in the relevant cohomology theory in function field arithmetic \cite[Corollary 15.5]
{Boeckle}. 

For Drinfeld cuspforms of rank three, B\"{o}ckle, Gr\"{a}f and Kaiser formulate related questions about slopes \cite[\S6]{BGK}.

We outline the main ideas of the proof of Theorem \ref{ThmMainIntro}. 
By the constancy theorem \cite[Theorem 1.1 (1)]{Ha_DMDia}, we may assume $r=1$.
In this case, 
by the formula of Bandini--Valentino \cite[(4.1)]{BV0},
a representing matrix of $U_t$ is written as
\[
U=B\diag(t^{a_0},\ldots,t^{a_{k-2}})
\]
with some $B\in M_{k-1}(\bF_q)$ and nonnegative integers $a_0<\cdots <a_{k-2}$.

The key input is Theorem \ref{ThmSymSlopes}, which shows that the $t$-adic valuation of any nonzero eigenvalue of a matrix of this type
is $\frac{1}{2}(a_i+a_j)$ for some $i,j$, provided that $B$ has a symmetric rank profile matrix. 
This can be viewed as a variant of \cite[Lemma 5]{BC} which Buzzard and Calegari used in order to compute the slopes of the $U_2$ operator 
acting on the space of 
$2$-adic overconvergent elliptic cuspforms of tame level one and weight zero.

We can write $B=L'RU'$ for some lower triangular invertible matrix $L'$, subpermutation matrix $R$ and upper unitriangular matrix $U'$. 
The matrix $R$ is called the rank profile matrix of $B$ \cite[Theorem 3]{DPS}. To prove Theorem \ref{ThmSymSlopes}, we pass to a matrix $\tilde{C}$
similar to $U$ and study the $t$-adic Newton polygon of its reciprocal characteristic polynomial
\[
P_{\tilde{C}}(x)=\det(\id-x\tilde{C}). 
\]
The coefficient of $x^s$ in this polynomial is, up to a sign, the sum of $s\times s$ principal minors of $\tilde{C}$. 
In general, we only have a lower bound for the $t$-adic valuation of such a principal minor. On the other hand, if $R$ is symmetric, then we can compute the exact 
$t$-adic 
valuations of some principal minors of $\tilde{C}$, which will give the vertices of the Newton polygon of $P_{\tilde{C}}(x)$.

Let $V$ be the $F$-subspace of the polynomial ring $F[X,Y]$ consisting of homogeneous polynomials of degree one, on which $g\in \mathit{GL}_2(F)$ acts by
\[
g\circ (X,Y)=(X,Y)g.
\]
For any $b\in \bF_q$ and $c\in \bF_q^\times$, let
\[
u_b=\begin{pmatrix}
	1&b\\0&1
\end{pmatrix},\quad v_c=\begin{pmatrix}
	1&-c\\c^{-1}&0
\end{pmatrix},\quad \tau=\begin{pmatrix}
1& 0\\ 0&t
\end{pmatrix}.
\]
Define elements $\cB,\cP\in F[\mathit{GL}_2(F)]$ by
\[
\cB=\sum_{b\in \bF_q}[u_b]-\sum_{c\in \bF_q^\times}[v_c],\quad \cP=1-\sum_{b\in \bF_q}[u_b].
\]

Next we will show that
the representing matrix $U$ is equal to a representing matrix of the action of $\cB\tau\in F[\mathit{GL}_2(F)]$
on the dual $(\Delta_{k-2}^\star)^\vee:=(\Sym^{k-2}(V^\vee))^\vee$ of the symmetric power representation of the dual of $V$ over $F$. 
Thus we reduce ourselves to computing $t$-adic valuations of eigenvalues of the action of $\cB\tau$ on each Jordan--H\"{o}lder factor of it.

By \cite[Corollary 2.4]{BGP}, any Jordan--H\"{o}lder factor of $(\Delta_{k-2}^\star)^\vee$ is a determinant twist of the tensor product of Frobenius twists 
of 
representations $(\Sym^m(V^\vee))^\vee$ for $m\in \{0,\ldots,p-1\}$. 
Let $(X^\vee, Y^\vee)$ be the dual basis of $V^\vee$ for $(X,Y)$.
Since $m<p$, we have the natural perfect pairing
\[
\begin{aligned}
\Sym^m(V)\times \Sym^m(V^\vee)\to F
\end{aligned}
\]
defined by $(X^{m-i}Y^i,(X^\vee)^{m-j}(Y^\vee)^j)\mapsto \delta_{i,j}(m-i)!i!$.
This induces a $\mathit{GL}_2(F)$-equivariant isomorphism
\[
(\Sym^m(V^\vee))^\vee\to \Sym^m(V)
\]
which is represented by an invertible diagonal matrix. 
Combining this with the key equality ${}^t\! (\cB\cP)=\cB\cP$ (Proposition \ref{PropKeySymmetry}), we can prove the required symmetry of the rank profile matrix 
for 
a representing matrix of the action of $\cB$ on each Jordan--H\"{o}lder factor 
(Corollary \ref{CorBSymRankProfile}).

In this way, we can compute slopes of $U_t$. More precisely, since $\cB$, $\cP$ and $\tau$ commute with $\begin{pmatrix}
	1&0\\0&d
\end{pmatrix}$ for any $d\in \bF_q^\times$, we can further reduce the problem to computing the finite slopes of matrices of the form
\[
B'\diag(t^{a'_0},\ldots,t^{a'_{l}}),
\]
where $B'$ has a symmetric rank profile matrix and the $a'_i$ are all congruent modulo $q-1$. 
Thus Theorem \ref{ThmSymSlopes} shows that the $t$-adic valuation of any nonzero eigenvalue of such a matrix is of the form
\[
N_1+N_2\left(\frac{q-1}{2}\right),\quad N_1,N_2\in \bZ.
\]
Now Theorem \ref{ThmMainIntro} (\ref{ThmMainIntro_HI})
follows since $\frac{q-1}{2}$ belongs to $\frac{1}{2}\bZ$, and it is an integer for odd $q$. 
On the other hand, \cite[Corollary 2.4]{BGP} provides control over symmetric powers and determinant twists appearing as factors of $(\Sym^{k-2}(V^\vee))^\vee$,
from which Theorem \ref{ThmMainIntro} (\ref{ThmMainIntro_Bound}) follows.

\subsection*{AI disclosure}
The initial proof strategy for Theorem \ref{ThmMainIntro} (\ref{ThmMainIntro_HI}) was developed through mathematical discussions with ChatGPT (OpenAI), 
using the GPT-6 Astra Pro model. The author 
checked the arguments, resolved ambiguities and 
wrote the initial manuscript. The AI model was subsequently used for proofreading and improving the exposition. The author is responsible for the mathematical 
content 
and the final text.

\subsection*{Acknowledgments}
This work was supported by JSPS KAKENHI Grant Number 26K06733.



\section{Symmetric rank profile matrices and half-integrality of slopes}\label{SecSymSlope}

Let $K$ be a field and let $t$ be an indeterminate. Let $d\geq 1$ be an integer. 
We adopt the convention that indexing of rows and columns in matrices starts with zero.
For a matrix $C\in M_d(K)$ and integers $u,v\in\{0,\ldots,d-1\}$, we denote by $C_{u,v}$ the $(u,v)$-entry of $C$, and let
\[
C_{\{0,\ldots,u\},\{0,\ldots,v\}}=(C_{i,j})_{i\in \{0,\ldots,u\},j\in \{0,\ldots,v\}}
\]
denote
the $(u+1)\times (v+1)$ northwest submatrix 
of $C$ consisting of the first $u+1$ rows and the first $v+1$ columns. Let 
\[
r_C(u,v)=\rank(C_{\{0,\ldots,u\},\{0,\ldots,v\}}).
\]
Also, let $r_C(u,v)=0$ when $u\notin \{0,\ldots,d-1\}$ or $v\notin \{0,\ldots,d-1\}$.
We denote the support of $C$ by 
\[
\Supp(C)=\{(i,j)\in\{0,\ldots,d-1\}^2\mid C_{i,j}\neq 0\}.
\]

\begin{lem}\label{LemSupport}
	Let $L,C,U\in M_d(K)$ be matrices such that 
	\begin{itemize}
		\item $L$ is lower triangular and invertible, and
		\item $U$ is upper triangular and invertible.
	\end{itemize}
	Then, for any $u,v\in \{0,\ldots,d-1\}$, we have
	\[
	r_{C}(u,v)=r_{LCU}(u,v).
	\]
\end{lem}
\begin{proof}
	From the assumptions on $L$ and $U$, we see that
	\[
	(LCU)_{\{0,\ldots,u\},\{0,\ldots,v\}}=L_{\{0,\ldots,u\},\{0,\ldots,u\}}C_{\{0,\ldots,u\},\{0,\ldots,v\}}U_{\{0,\ldots,v\},\{0,\ldots,v\}}
	\]
	for any $u,v\in \{0,\ldots,d-1\}$, where $L_{\{0,\ldots,u\},\{0,\ldots,u\}}$ is lower triangular invertible and $U_{\{0,\ldots,v\},\{0,\ldots,v\}}$ is upper 
	triangular invertible. This yields the equality of the lemma. 
\end{proof}

\begin{dfn}\label{DfnQuasiMonomialMatrix}
	A matrix $C\in M_d(K)$ is said to be a subpermutation matrix if every entry of $C$ is zero or one, and $C$ has at most one nonzero entry in
	each row and column.
\end{dfn}

A key notion in this paper is the rank profile matrix \cite[\S2.1]{DPS}.

\begin{lem}\label{LemLMUDecomp}
	For any $C\in M_d(K)$, there exist $L,U,R\in M_d(K)$ such that $L$ is lower triangular and invertible, $U$ is upper unitriangular, $R$ is a subpermutation 
	matrix and 
	$C=LRU$. The matrix $R$ is called the rank profile matrix of $C$.
\end{lem}
\begin{proof}
	We may assume $C\neq O$. Let $i\in \{0,\ldots,d-1\}$ be the smallest integer such that the $i$-th row of $C$ is nonzero, and let 
	$j\in 
	\{0,\ldots, d-1\}$ be the smallest integer such that $C_{i,j}$ is nonzero. Then $C_{i,l}=0$ for any $l<j$ and $C_{l,j}=0$ for any $l<i$. By multiplying $C$ 
	from 
	the left by a lower
	unitriangular matrix and from the right by an upper unitriangular matrix, we can transform $C$ to satisfy $C_{i,l}=0$ for any $l\neq j$ and 
	$C_{l,j}=0$ for 
	any $l\neq i$. By multiplying an invertible diagonal matrix from the left, we may assume $C_{i,j}=1$.
	Repeating this procedure on the remaining rows and columns proves the lemma.
\end{proof}

\begin{lem}\label{LemSuppVSRankFunction}
	Let $R\in M_d(K)$ be a subpermutation matrix. Then, for any $i,j\in \{0,\ldots,d-1\}$, we have
	\begin{equation}\label{EqnSuppRankFunction}
		\begin{aligned}
	(i,j)&\in \Supp(R) \\
	&\Leftrightarrow r_R(i,j)-r_R(i-1,j)-r_R(i,j-1)+r_R(i-1,j-1)=1.
	\end{aligned}
	\end{equation}
	In particular, for any $C\in M_d(K)$, its rank profile matrix is unique.
\end{lem}
\begin{proof}
	Set $R_{i,j}=0$ when $i$ or $j$ is not in $\{0,\ldots, d-1\}$.
	For any $u,v\in \bZ$, define $\vep_{u,v}\in \bZ$ by
	\[
	\vep_{u,v}=\left\{
	\begin{array}{ll}
		1 & (R_{u,v}\neq 0),\\
		0 & (R_{u,v}=0).
	\end{array}
	\right.
	\]
	Note that $\vep_{u,v}\in \bZ$, while $R_{u,v}\in K$.
	From the definition of a subpermutation matrix, we see that for any $i,j\in\{0,\ldots,d-1\}$, we have
	\[
	\begin{aligned}
	r_R(i,j)&=|\Supp(R)\cap (\{0,\ldots,i\}\times \{0,\ldots,j\})|=\sum_{u=0}^i\sum_{v=0}^j \vep_{u,v}.
	\end{aligned}
	\]
	This yields
	\[
	\begin{aligned}
		r_R(i,j)-r_R(i-1,j)&=\sum_{v=0}^j \vep_{i,v},\\
		r_R(i,j-1)-r_R(i-1,j-1)&=\sum_{v=0}^{j-1} \vep_{i,v}.
	\end{aligned}
	\]
	Hence we obtain
	\[
	 r_R(i,j)-r_R(i-1,j)-r_R(i,j-1)+r_R(i-1,j-1)=\vep_{i,j},
	\]
	from which the equivalence (\ref{EqnSuppRankFunction}) follows. If $C=LRU$ is a decomposition as in Lemma \ref{LemLMUDecomp}, then Lemma \ref{LemSupport}
	yields
	$r_R(u,v)=r_C(u,v)$.
	Since a subpermutation matrix is determined by its support, (\ref{EqnSuppRankFunction}) gives
	the last assertion.
\end{proof}

\begin{lem}\label{LemSymmRankProfile}
	Let $L,B,U\in M_d(K)$ be matrices satisfying the following conditions.
	\begin{itemize}
		\item $L$ is lower triangular and invertible.
		\item $U$ is upper triangular and invertible.
		\item $LBU$ is symmetric.
	\end{itemize}
	Then the rank profile matrix $R$ of $B$ is symmetric.
\end{lem}
\begin{proof}
	The assumption on the symmetry of $LBU$ and Lemma 
	\ref{LemSupport} show
	\begin{equation*}
		r_R(u,v)=r_{LBU}(u,v)=r_{LBU}(v,u)=r_R(v,u)
	\end{equation*}
	for any $u,v\in \{0,\ldots,d-1\}$. Since $R$ is a subpermutation matrix, 
	this symmetry of $r_R(u,v)$ and Lemma \ref{LemSuppVSRankFunction} show that the set $\Supp(R)$ is stable under the involution $(i,j)\mapsto (j,i)$.
	This means that $R$ is symmetric.
\end{proof}

We fix an algebraic closure $\overline{K(t)}$ of $K(t)$, an algebraic closure $\overline{K((t))}$ of $K((t))$ and an embedding $\iota:\overline{K(t)}\to 
\overline{K((t))}$ of $K(t)$-algebras. 
Let $v_t$ be the normalized additive $t$-adic valuation on $\overline{K((t))}$ satisfying $v_t(t)=1$. For any element $\alpha\in \overline{K(t)}$, 
we define the $t$-adic valuation of $\alpha$ by $v_t(\iota(\alpha))$.
For any integer $d\geq 1$, any matrix $C\in M_d(K(t))$ and any eigenvalue $\alpha\in \overline{K(t)}$ of $C$, we refer to $v_t(\iota(\alpha))$ as a slope of $C$.

The following theorem will enable us to compute the slopes of $U_t$ acting on the space of Drinfeld cuspforms.

\begin{thm}\label{ThmSymSlopes}
	 Let $a_0, a_1,\cdots,a_{d-1}$ be nonnegative integers satisfying
	\[
	a_0<a_1<\cdots<a_{d-1}.
	\] 
	Let $D=\diag(t^{a_0},\ldots,t^{a_{d-1}})$ and $B\in M_d(K)$.
	Suppose that the rank profile matrix $R$ of $B$ is symmetric.
	Then, every finite slope of the matrix $BD$ equals
	\[
	\frac{1}{2}(a_i+a_j)
	\] 
	for some $i,j\in \{0,\ldots,d-1\}$.
\end{thm}
\begin{proof}
Note that, since $BD\in M_d(K[t])$, the multiset of its slopes is independent of the choice of an embedding $\iota:\overline{K(t)}
\to\overline{K((t))}$: it can be read off from the $t$-adic Newton polygon of the reciprocal characteristic polynomial
\[
\det(\id-xBD)
\] 	
in an indeterminate $x$.

We can write $B=L'RU'$ in $M_d(K)$ with some lower triangular invertible matrix $L'$ and upper 
unitriangular matrix 
$U'$.
If $R=O$, then $B=O$ and the theorem is immediate.
Thus we may assume $R\neq O$.
Let 
\[
I^{\times}=\{i\in \{0,\ldots,d-1\}\mid \text{the $i$-th row of $R$ is nonzero}\}\neq \emptyset.
\]
For any $i\in I^{\times}$, there exists a unique $\sigma(i)\in \{0,\ldots,d-1\}$ satisfying $(i,\sigma(i))\in \Supp(R)$. 
Then the symmetry of $R$ implies $(\sigma(i),i)\in \Supp(R)$, which means $\sigma(i)\in I^{\times}$ and $\sigma(\sigma(i))=i$. Thus the map $\sigma$ defines 
an involution $\sigma:I^{\times}\to I^{\times}$.

Let $t^{1/2}$ be a square root of $t$ in $\overline{K((t))}$ and consider the natural embeddings 
\[
K(t)\to K((t))\to K((t^{1/2}))
\]
of $K$-subalgebras of $\overline{K((t))}$.
Let
\[
D^{1/2}=\diag(t^{a_0/2},\ldots,t^{a_{d-1}/2})\in M_d(K[[t^{1/2}]])
\]
so that $(D^{1/2})^2=D$ and $(D^{1/2})^{-1}=\diag(t^{-a_0/2},\ldots,t^{-a_{d-1}/2})$. Let
\[
\tilde{L}'=D^{1/2}L'(D^{1/2})^{-1},\quad \tilde{R}=D^{1/2}RD^{1/2},\quad \tilde{U}'=(D^{1/2})^{-1}U'D^{1/2}.
\]
Since $BD$ and $(\tilde{L}')^{-1}D^{1/2}(BD)(D^{1/2})^{-1}\tilde{L}'=\tilde{R}\tilde{U}'\tilde{L}'$ have the same reciprocal characteristic polynomial, it is 
enough to compute the $t$-adic 
valuation of
each eigenvalue of the matrix $\tilde{C}:=\tilde{R}\tilde{U}'\tilde{L}'$.

Since the $a_i$ are strictly increasing, we have
\begin{equation}\label{EqnUtLtPrincipalUnit}
\tilde{U}'\in \id+t^{1/2}M_d(K[[t^{1/2}]]),\quad \tilde{L}'\in L'_0+t^{1/2}M_d(K[[t^{1/2}]]),
\end{equation}
where $L'_0=\diag(L'_{0,0},\ldots, L'_{d-1,d-1})$. Note that $L'_0\in \mathit{GL}_d(K)$ is an invertible diagonal matrix.
Moreover, the matrix $\tilde{R}$ satisfies $\Supp(\tilde{R})=\Supp(R)$. For any $i\in I^{\times}$, we have $(i,\sigma(i))\in \Supp(\tilde{R})$. 
Let 
\[
b_i:=\frac{a_i+a_{\sigma(i)}}{2}\quad (i\in I^\times).
\]
Then we can write
\begin{equation*}
	\tilde{R}_{i,j}=\left\{
	\begin{array}{ll}
		t^{b_i} & (i\in I^\times,\ j=\sigma(i)),\\
		0 & (\text{otherwise}).
	\end{array}
	\right.
\end{equation*}

Now we show that the multiset of $t$-adic valuations of nonzero eigenvalues of $\tilde{C}$ is $\{b_i\mid i\in I^\times\}$.
Let $\beta_1<\cdots <\beta_m$ be the distinct values in this multiset.
Decompose $I^\times$ as
\[
I^\times=\coprod_{l=1}^m I^\times_l,\quad I^\times_l=\{i\in I^\times\mid b_i=\beta_l\}.
\]
Let $n_l=|I^\times_l|$.

Let
\[
P_{\tilde{C}}(x):=\det(\id-x\tilde{C})
\]
be the reciprocal characteristic polynomial of $\tilde{C}$.
For any $I\subseteq \{0,\ldots,d-1\}$, 
let $\tilde{C}_{I}$ be the principal submatrix of $\tilde{C}$ for $I$. 
Recall that, for any $s\in \{1,\ldots, d\}$, the coefficient of $x^s$ in $P_{\tilde{C}}(x)$ 
equals
\[
(-1)^s \sum_{|I|=s} \det(\tilde{C}_{I}).
\]

Take any $s\in \{1,\ldots,d\}$ and $I\subseteq \{0,\ldots,d-1\}$ of cardinality $s$. Write
\[
\tilde{U}'\tilde{L}'=(e_{i,j})_{i,j},\quad e_{i,j}\in K[[t^{1/2}]].
\] 
Then 
\begin{equation*}
	\tilde{C}_{i,j}=\left\{
	\begin{array}{ll}
		t^{b_i} e_{\sigma(i),j} & (i\in I^\times),\\
		0 & (i\notin I^\times).
	\end{array}
	\right.
\end{equation*}
This implies 
\begin{equation}\label{EqnValDetCt}
	\left\{
	\begin{array}{ll}
		\det(\tilde{C}_I)=0 & (I\not\subseteq I^\times),\\
v_t(\det(\tilde{C}_{I}))\geq \sum_{i\in I}b_i & (I\subseteq I^\times).
\end{array}
\right.
\end{equation}
In particular, if $s\leq |I^\times|$, then the $t$-adic valuation of the coefficient of $x^s$ in $P_{\tilde{C}}(x)$ is no less than the sum of the $s$ 
smallest values in the 
multiset $\{b_i\mid i\in I^\times\}$. If $s> |I^\times|$, then that coefficient is zero.

Suppose that
\[
s=n_1+\cdots+n_l,\quad I=I^\times_1\sqcup\cdots\sqcup I^\times_l
\] 
for some $l\in\{1,\ldots,m\}$, so that $I\subseteq I^\times$. 
For any subset $I'\subseteq \{0,\ldots,d-1\}$ of cardinality $s$
other than $I$, (\ref{EqnValDetCt}) yields
\begin{equation}\label{EqnNPNonVertex}
	\left\{
		\begin{array}{ll}
		\det(\tilde{C}_{I'})=0 & (I'\not\subseteq I^\times),\\
v_t(\det(\tilde{C}_{I'}))\geq \sum_{i\in I'}b_i>\sum_{i\in I}b_i & (I'\subseteq I^\times).
\end{array}
\right.
\end{equation}
On the other hand, since $\sigma$ is an involution, we have $b_{\sigma(i)}=b_i$ and our choice of $I$ implies $\sigma(I)=I$. Since
\[
\det(\tilde{C}_{I})=\det((t^{b_i}e_{\sigma(i),j})_{i,j\in I})=\det((t^{b_{\sigma(i)}}e_{\sigma(i),j})_{i,j\in I})
\] 
and the matrix in the last determinant is the one obtained from $(t^{b_{i}}e_{i,j})_{i,j\in I}$ by a row permutation, we have
\[
\det(\tilde{C}_{I})=\pm\det((t^{b_{i}}e_{i,j})_{i,j\in I})=\pm t^{\sum_{i\in I}b_i} \det((e_{i,j})_{i,j\in I}).
\]
By (\ref{EqnUtLtPrincipalUnit}), we have $\tilde{U}'\tilde{L}'\in L'_0+t^{1/2}M_d(K[[t^{1/2}]])$ and thus 
\[
(e_{i,j})_{i,j\in I}\in (L'_0)_I+t^{1/2}M_s(K[[t^{1/2}]]), 
\]
where $(L'_0)_I$ is the principal submatrix of $L'_0$ for $I$.
Since $(L'_0)_I$ is also an invertible diagonal matrix in $\mathit{GL}_s(K)$, we have 
\[
\det((e_{i,j})_{i,j\in I})\in K[[t^{1/2}]]^\times.
\]
Hence we obtain
\begin{equation}\label{EqnNPVertex}
v_t(\det(\tilde{C}_{I}))=\sum_{i\in I}b_i.
\end{equation}
Then (\ref{EqnNPNonVertex}) and (\ref{EqnNPVertex}) imply that, if $s=n_1+\cdots+n_l$ for some $l\in\{1,\ldots,m\}$, then 
the $t$-adic valuation of the coefficient of $x^s$ in $P_{\tilde{C}}(x)$ is equal to 
the sum of the $s$ smallest values in the multiset $\{b_i\mid i\in I^\times\}$. Therefore, the $t$-adic Newton polygon of $P_{\tilde{C}}(x)$ has slope 
multiset
$\{b_i\mid i\in I^\times\}$. This concludes the proof of the theorem.
\end{proof}



\section{Jordan--H\"{o}lder factors of symmetric power representations of $\mathit{GL}_2(F)$}\label{SecJH}

Let $F[\mathit{GL}_2(F)]$ be the group ring of $\mathit{GL}_2(F)$ over $F$. For any $g\in \mathit{GL}_2(F)$, we denote by $[g]\in F[\mathit{GL}_2(F)]$
the corresponding element. We denote $[g]$ also by $g$ when no confusion can arise. For any representation
\[
\rho:\mathit{GL}_2(F)\to \mathit{GL}_d(F),
\]
we denote by the same symbol $\rho$ the $F$-algebra homomorphism
\[
F[\mathit{GL}_2(F)]\to M_{d}(F)
\]
induced by $\rho$. Similarly, we denote by
\[
F[\mathit{GL}_2(F)]\to F[\mathit{GL}_2(F)],\quad \cG\mapsto {}^t\!\cG
\]
the $F$-linear map induced by the transpose 
\[
\mathit{GL}_2(F)\to \mathit{GL}_2(F),\quad g\mapsto {}^t\!g.
\]
Then we have ${}^t\!(\cG\cH)={}^t\!\cH {}^t\!\cG$ for any $\cG,\cH\in F[\mathit{GL}_2(F)]$.

For any integer $m\geq 0$, let $\Delta_m$ be the $F$-subspace of the polynomial ring $F[X,Y]$ consisting of homogeneous polynomials of degree $m$.
As in \cite{Boeckle, Bon}, we consider the left action of $\mathit{GL}_2(F)$ on $\Delta_m$ defined by
\[
g\circ P(X,Y):= P((X,Y)g)
\]
for any $g\in \mathit{GL}_2(F)$ and $P(X,Y)\in \Delta_m$.
We denote by $r_m(g)$ the representing matrix of $g$ for the action $\circ$ with respect to the basis
\[
\Sigma_m=(X^m, X^{m-1}Y,\ldots,XY^{m-1},Y^m)
\]
of $\Delta_m$ so that
\[
\begin{aligned}
	g\circ (X^m, X^{m-1}Y,\ldots,XY^{m-1},Y^m)=(X^m, X^{m-1}Y,\ldots,XY^{m-1},Y^m)r_m(g).
\end{aligned}
\]

On the other hand, we denote by $\Delta^\star_m$ the $F$-vector space $\Delta_m$ equipped with the left action of $\mathit{GL}_2(F)$
defined by
\[
g\star P(X,Y):= P((X,Y){}^t\! g^{-1})
\]
for any $g\in \mathit{GL}_2(F)$ and $P(X,Y)\in \Delta_m$. This is the normalization adopted in \cite[\S2.1]{BGP}.
We denote by $r^\star_m(g)$ the representing matrix of $g$ for the action $\star$ 
with respect to the basis $\Sigma_m$.
Using the automorphism
\[
\varphi_d:\mathit{GL}_d(F)\to \mathit{GL}_d(F),\quad g\mapsto {}^t\!g^{-1},
\]
we can write
\[
r^\star_m=r_m\circ \varphi_2:\mathit{GL}_2(F)\to \mathit{GL}_{m+1}(F).
\]

Let 
\[
\Delta_m^\vee:=\Hom_F(\Delta_m,F),\quad (\Delta_m^\star)^\vee:=\Hom_F(\Delta^\star_m,F)
\]
on which we define left actions of $\mathit{GL}_2(F)$ by
\begin{equation}\label{EqnActionDualStar}
\begin{aligned}
(g\circ\omega)(P(X,Y))&:=\omega(g^{-1}\circ P(X,Y))=\omega(P((X,Y)g^{-1})),\\
(g\star\omega^\star)(P(X,Y))&:=\omega^\star(g^{-1}\star P(X,Y))=\omega^\star(P((X,Y){}^t\! g))
\end{aligned}
\end{equation}
for any $\omega\in \Delta_m^\vee$ and $\omega^\star\in (\Delta_m^\star)^\vee$.

For any $j\in \{0,\ldots,m\}$, let $(X^{m-j}Y^j)^\vee$ be the unique element of $\Delta_m^\vee$ and $(\Delta_m^\star)^\vee$ satisfying
\[
(X^{m-j}Y^j)^\vee(X^{m-i}Y^i)=\delta_{i,j}
\]
for any $i\in \{0,\ldots,m\}$. We denote by $\rho_m(g)$ the representing matrix of the action $\circ$ on $\Delta_m^\vee$ with respect to the basis
\[
\Sigma^\vee_m=((X^m)^\vee, (X^{m-1}Y)^\vee,\ldots,(XY^{m-1})^\vee,(Y^m)^\vee)
\]
of $\Delta_m^\vee$ so that
\[
\begin{aligned}
	g\circ &((X^m)^\vee, (X^{m-1}Y)^\vee,\ldots,(XY^{m-1})^\vee,(Y^m)^\vee)\\&=((X^m)^\vee, (X^{m-1}Y)^\vee,\ldots,(XY^{m-1})^\vee,(Y^m)^\vee)\rho_m(g).
\end{aligned}
\]
We define $\rho_m^\star(g)$ similarly for $(\Delta_m^\star)^\vee$.
Then they give homomorphisms
\[
\rho_m:\mathit{GL}_2(F)\to \mathit{GL}_{m+1}(F),\quad \rho^\star_m:\mathit{GL}_2(F)\to \mathit{GL}_{m+1}(F)
\]
satisfying 
\begin{equation*}
	\rho_m=\varphi_{m+1}\circ r_m,\quad \rho^\star_m=\varphi_{m+1}\circ r^\star_m,\quad \rho_m^\star=\rho_m\circ \varphi_2.
	\end{equation*}
Note that we have isomorphisms of left $F[\mathit{GL}_2(F)]$-modules
\[
\Delta_1^\star\to (\Delta_1)^\vee,\quad \Delta_m^\star\to \Sym^m((\Delta_1)^\vee).
\]
However, since the field $F$ is of characteristic $p$, duality does not in general commute with $m$-th symmetric powers.
As we will see in the following lemma, this commutativity holds for $m<p$.

For any $m\in \{0,\ldots,p-1\}$, let
\[
H_m=\diag\left(\binom{m}{0},\ldots, \binom{m}{m}\right)\in \mathit{GL}_{m+1}(\bF_p).
\]

\begin{lem}\label{LemIsomStarCircDual}
	For any $m\in \{0,\ldots,p-1\}$, we have an isomorphism of left $F[\mathit{GL}_2(F)]$-modules
	\[
	h_m:(\Delta_m^\star)^\vee\to \Delta_m,\quad  (X^{m-i}Y^i)^\vee\mapsto \binom{m}{i} X^{m-i}Y^i.
	\]
\end{lem}
\begin{proof}
	Write $V=\Delta_1$ and $V^\vee=\Delta_1^\vee$. Then we have an isomorphism of left $F[\mathit{GL}_2(F)]$-modules
	\[
	\Phi: \Delta^\star_m\to \Sym^m(V^\vee),\quad X^{m-i}Y^i\mapsto (X^\vee)^{m-i}(Y^\vee)^i.
	\]
	On the other hand, since $m<p$, we have the perfect pairing
	\[
	\langle-,-\rangle:\Sym^m(V)\times \Sym^m(V^\vee)\to F
	\]
	defined by
	\[
	\langle v_1\cdots v_m,l_1\cdots l_m\rangle=\sum_{\sigma\in S_m} \prod_{i=1}^m l_i(v_{\sigma(i)}),
	\]
	where $v_i\in V$, $l_i\in V^\vee$ and $S_m$ is the symmetric group of degree $m$. For any $g\in \mathit{GL}_2(F)$, we have
	\[
	\langle g\circ v_1\cdots v_m,l_1\cdots l_m\rangle=\langle v_1\cdots v_m,g^{-1}\circ l_1\cdots l_m \rangle.
	\]
	This induces an isomorphism of left $F[\mathit{GL}_2(F)]$-modules
	\[
	\begin{aligned}
	\Psi:\Sym^m(V)&\to (\Sym^m(V^\vee))^\vee,\\
	v_1\cdots v_m&\mapsto (l_1\cdots l_m\mapsto \langle v_1\cdots v_m,l_1\cdots l_m \rangle).
	\end{aligned}
	\]
	Since 
	\[
	\langle X^{m-i}Y^i,(X^\vee)^{m-j}(Y^\vee)^j\rangle=\delta_{i,j}(m-i)!i!,
	\]
	we obtain an isomorphism of left $F[\mathit{GL}_2(F)]$-modules
	\[
	\Phi^\vee\circ \Psi: \Sym^m(V)\to (\Delta_m^\star)^\vee,\quad X^{m-i}Y^i\mapsto (m-i)!i!(X^{m-i}Y^i)^\vee.
	\]
	Since $m!\in F^\times$, the map $h_m=m!(\Phi^\vee\circ \Psi)^{-1}$ is also an isomorphism of left $F[\mathit{GL}_2(F)]$-modules.
\end{proof}

For any integers $s\geq 0$ and $d\geq 1$, the $p^s$-th power Frobenius endomorphism $\Fr^{(s)}:F\to F$ defines the homomorphism
\[
\Fr^{(s)}:\mathit{GL}_d(F)\to \mathit{GL}_d(F),\quad (c_{i,j})_{i,j}\mapsto (c_{i,j}^{p^s})_{i,j}.
\] 
We denote by $\Delta_m^{(s)}$ the $F$-vector space $\Delta_m$ on which
$\mathit{GL}_2(F)$ acts by the composition 
\[
\Fr^{(s)}\circ r_m: \mathit{GL}_2(F)\to \mathit{GL}_{m+1}(F)\to \mathit{GL}_{m+1}(F)
\] 
with respect to the basis $\Sigma_m$. 
We define $(\Delta_m^\star)^{(s)}$ similarly.
We also define $((\Delta_m^\star)^\vee)^{(s)}$ using the basis $\Sigma_m^\vee$.

Let $m\geq 0$ and $a$ be integers. Let
\[
m=m_0+m_1 p+\cdots+ m_\mu p^\mu,\quad m_\nu\in\{0,\ldots,p-1\}
\]
be the base-$p$ expansion of $m$. Put
\[
L_{m,a}:=\mathrm{det}^{a}\otimes \bigotimes_{\nu=0}^\mu \Delta_{m_\nu}^{(\nu)},\quad 
L_{m,a}^\star:=\mathrm{det}^{a}\otimes \bigotimes_{\nu=0}^\mu (\Delta_{m_\nu}^\star)^{(\nu)}.
\]
We also put
\[
W_{m,a}^\star:=\mathrm{det}^{a}\otimes \bigotimes_{\nu=0}^\mu ((\Delta_{m_\nu}^\star)^\vee)^{(\nu)}.
\]
They are $F$-vector spaces of dimension $d_m:=\prod_{\nu=0}^\mu(m_\nu+1)$.

For any $\underline{i}=(i_\nu)_{\nu\in \{0,\ldots,\mu\}}\in \{0,\ldots,m_0\}\times \cdots\times \{0,\ldots,m_\mu\}$, let
\[
\begin{aligned}
f_{\underline{i}}&:=X^{m_0-i_0}Y^{i_0}\otimes \cdots \otimes X^{m_\mu-i_\mu}Y^{i_\mu},\\
f_{\underline{i}}^\wedge&:=(X^{m_0-i_0}Y^{i_0})^\vee\otimes \cdots \otimes (X^{m_\mu-i_\mu}Y^{i_\mu})^\vee.
\end{aligned}
\]
We denote by $T_m$ the basis of $L_{m,a}$ and $L_{m,a}^\star$ defined by arranging
\[
f_{\underline{i}},\quad \underline{i}\in \{0,\ldots,m_0\}\times \cdots \times\{0,\ldots,m_\mu\}
\]
in the colexicographic order: for $\underline{i}=(i_\nu)_\nu$ and $\underline{i}'=(i'_\nu)_\nu$, we say $\underline{i}<\underline{i}'$ if there exists $\nu_0$ such 
that $i_{\nu_0}<i'_{\nu_0}$ and for any $\nu>\nu_0$ we have 
$i_\nu=i'_\nu$. Put
\[
w(\underline{i})=\sum_{\nu=0}^\mu i_\nu p^\nu,
\]
so that
\begin{equation}\label{EqnRangeOfWeight}
	w(\underline{i})\in \bZ,\quad 0\leq w(\underline{i})\leq m.
\end{equation}
Since $0\leq i_\nu\leq p-1$, we have
\begin{equation}\label{EqnCoLexVSWeight}
\underline{i}<\underline{i}' \quad\Leftrightarrow \quad w(\underline{i})<w(\underline{i}').
\end{equation}

Similarly, we denote by $T^\wedge_m$ the basis of $W_{m,a}^\star$ defined by arranging
\[
f_{\underline{i}}^\wedge,\quad \underline{i}\in \{0,\ldots,m_0\}\times \cdots \times\{0,\ldots,m_\mu\}
\]
in the colexicographic order, and by $T^\vee_m$ the dual basis of $(L_{m,a}^\star)^\vee$ for $T_m$ consisting of 
\[
f_{\underline{i}}^\vee,\quad \underline{i}\in \{0,\ldots,m_0\}\times \cdots \times\{0,\ldots,m_\mu\}
\]
in the colexicographic order.
Let
\[
r_{m,a},r^\star_{m,a},\rho^\star_{m,a}:\mathit{GL}_2(F)\to \mathit{GL}_{d_m}(F)
\]
be the matrix representations of the actions of $\mathit{GL}_2(F)$ on $L_{m,a}$, $L_{m,a}^\star$, $W_{m,a}^\star$
with respect to the bases $T_m$, $T_m$, $T_m^\wedge$, respectively.

\begin{lem}\label{LemLstarDual}
	There exists an isomorphism of left $F[\mathit{GL}_2(F)]$-modules
	\[
	(L_{m,a}^\star)^\vee\to W_{m,-a}^\star,\quad f_{\underline{i}}^\vee\mapsto f_{\underline{i}}^\wedge.
	\]
\end{lem}
\begin{proof}
	For any $C\in \mathit{GL}_{m_\nu+1}(F)$, we have $\Fr^{(s)}({}^t\! C^{-1})={}^t\!(\Fr^{(s)}(C))^{-1}$. 
	Thus we have an equality of left $F[\mathit{GL}_2(F)]$-modules
	\[
	\mathrm{det}^{-a}\otimes \bigotimes_{\nu}((\Delta_{m_\nu}^\star)^{(\nu)})^\vee = \mathrm{det}^{-a}\otimes\bigotimes_{\nu}((\Delta_{m_\nu}
	^\star)^\vee)^{(\nu)}.
	\]
	On the other hand, since duality commutes with tensor products, we have an isomorphism of left $F[\mathit{GL}_2(F)]$-modules
	\[
	(\mathrm{det}^{a}\otimes \bigotimes_{\nu}((\Delta_{m_\nu}^\star)^{(\nu)}))^\vee \to \mathrm{det}^{-a}\otimes \bigotimes_{\nu}((\Delta_{m_\nu}
	^\star)^{(\nu)})^\vee,\quad f_{\underline{i}}^\vee\mapsto f_{\underline{i}}^\wedge.
	\]
	Hence the lemma follows.
\end{proof}

By representing the action of $g\in \mathit{GL}_2(F)$ by the bases $T_{m}^\vee$ of $(L_{m,a}^\star)^\vee$ and $T_m^\wedge$ of $W_{m,-a}^\star$,
the definitions and Lemma \ref{LemLstarDual} yield
\begin{equation}\label{EqnSymmetricPhiTwists}
	r^\star_{m,-a}=r_{m,a}\circ \varphi_2,\quad \rho^\star_{m,-a}=\varphi_{d_m}\circ r^\star_{m,a}.
\end{equation}

On the other hand, the tensor product of Frobenius twists of isomorphisms $h_{m_\nu}$ of Lemma \ref{LemIsomStarCircDual} gives
an isomorphism of left $F[\mathit{GL}_2(F)]$-modules
\[
h_{m}: W_{m,a}^\star\to L_{m,a}.
\]
Let $H_{m}$ be the representing matrix of $h_{m}$ with respect to the basis $T_{m}^\wedge$ of $W^\star_{m,a}$ and $T_m$ of $L_{m,a}$.
We have
\[
H_m=\diag_{\underline{i}}\left(\prod_{\nu=0}^\mu \binom{m_\nu}{i_\nu}\right),
\]
where the $\underline{i}$ are ordered colexicographically.
Note that $H_m$ is an invertible diagonal matrix.
For $m\in \{0,\ldots,p-1\}$, the matrix $H_m$ is the same as the one we defined in that case.
Since the map $h_m$ is $\mathit{GL}_2(F)$-equivariant, we obtain
\begin{equation}\label{EqnRmvsRhostarm}
r_{m,a}(g)H_{m}=H_{m}\rho^\star_{m,a}(g)
\end{equation}
for any $g\in \mathit{GL}_2(F)$.

\begin{lem}\label{LemSelfDualW}
	Let $m\in \bZ_{\geq 0}$ and $a\in \bZ$. Then we have
	\[
	{}^t\!(\rho_{m,a}^\star(g))H_{m}=H_{m} \rho_{m,a}^\star({}^t\!g),\quad {}^t\!(\rho_{m,a}^\star(\cG))H_{m}=H_{m} \rho_{m,a}^\star({}^t\!\cG)
	\]
	for any $g\in \mathit{GL}_2(F)$ and $\cG\in F[\mathit{GL}_2(F)]$. 
\end{lem}
\begin{proof}
	It is enough to show the first equality.
	By (\ref{EqnSymmetricPhiTwists}), we have
	\[
		{}^t\!(\rho_{m,a}^\star(g))H_{m}=r^\star_{m,-a}(g)^{-1}H_m=	(r_{m,a}({}^t\!g^{-1}))^{-1}H_{m}=r_{m,a}({}^t\!g)H_{m},
	\]
	which equals $H_{m}\rho^\star_{m,a}({}^t\!g)$ by (\ref{EqnRmvsRhostarm}).
\end{proof}

For any $d\in \bF_q^\times$, let
\[
y_d=\begin{pmatrix}
	1&0\\0&d
\end{pmatrix}.
\]
For any $e\in \{0,\ldots, q-2\}$, let
\[
\begin{gathered}
J_e:=\{\underline{i}\in \{0,\ldots,m_0\}\times \cdots \times\{0,\ldots,m_\mu\}\mid w(\underline{i})\equiv e\bmod q-1\},\\
W^\star_{m,a,e}:=\langle f^\wedge_{\underline{i}}\mid \underline{i}\in J_e\rangle.
\end{gathered}
\]
From (\ref{EqnActionDualStar}), we see that 
$W^\star_{m,a,e}$ is the maximal $F$-subspace of $W^\star_{m,a}$ such that $y_d$ acts by $d^{a+e}$ for any $d\in \bF_q^\times$.
Let 
\[
T_{m,e}^\wedge=(f_{\underline{i}}^\wedge)_{\underline{i}\in J_e}, 
\]
where we give $J_e$ the colexicographic order. Then $T_{m,e}^\wedge$ forms a basis of 
the $F$-vector space $W^\star_{m,a,e}$. 
Let $H_{m,J_e}$
be the principal submatrix of $H_m$ for $J_e$. Since $H_m$ is an invertible diagonal matrix, so is $H_{m,J_e}$.

Let $\cG\in F[\mathit{GL}_2(F)]$. Suppose that $\cG$ commutes with $y_d$ for any $d\in \bF_q^\times$.
Then the action of $\cG$ on $W^\star_{m,a}$ stabilizes $W^\star_{m,a,e}$ for any $e\in \{0,\ldots, q-2\}$. We denote by
$\rho_{m,a,e}^\star(\cG)$ the representing matrix of the action of $\cG$ on $W^\star_{m,a,e}$ with respect to its basis $T_{m,e}^\wedge$.
Then we have
\begin{equation}\label{EqnEBlockSubmatrix}
\rho^\star_{m,a,e}(\cG)=\rho^\star_{m,a}(\cG)_{J_e},
\end{equation}
where $\rho^\star_{m,a}(\cG)_{J_e}$ is the principal submatrix of $\rho^\star_{m,a}(\cG)$ for $J_e$.
Since ${}^t\!y_d=y_d$, we see that ${}^t\! \cG$ also commutes with the $y_d$ and thus $\rho_{m,a,e}^\star({}^t\!\cG)$ is well-defined.
By Lemma \ref{LemSelfDualW} and (\ref{EqnEBlockSubmatrix}), we obtain
\begin{equation}\label{EqnSeldual-e-part}
	{}^t\!(\rho_{m,a,e}^\star(\cG))H_{m,J_e}=H_{m,J_e} \rho_{m,a,e}^\star({}^t\!\cG).
\end{equation}

For any integer $n\geq 0$, let $n_0\in\{0,\ldots,p-1\}$ satisfy $n\equiv n_0\bmod p$. 
Define $\cE(n)\subseteq \bZ$ inductively by
\[
\cE(n)=\left\{
\begin{array}{ll}
	\{0\}& (n\in\{0,\ldots,p-1\}),\\
	p\cE\left(\frac{n-n_0}{p}\right) & (n\geq p\text{ and }n_0=p-1),\\
		p\cE\left(\frac{n-n_0}{p}\right)\sqcup \left(n_0+1+	p\cE\left(\frac{n-n_0-p}{p}\right)\right)\ & (n\geq p\text{ and }n_0\leq p-2).
	\end{array}
\right.
\]
An induction on $n$ shows that every $m\in n-2\cE(n)$ satisfies 
\begin{equation}\label{Eqnn-mhalfWellDef}
0\leq m\leq n,\quad m\equiv n\bmod 2.	
\end{equation}

\begin{prop}\label{PropJH}
	For any integer $n\geq 0$, the Jordan--H\"{o}lder factors of the left $F[\mathit{GL}_2(F)]$-module $(\Delta_n^\star)^\vee$
	are
	\[
	\left\{W^\star_{m,\frac{n-m}{2}}\ \middle|\ m\in n-2\cE(n)\right\}.
	\]
\end{prop}
\begin{proof}
	By \cite[Lemma 2.1 and Corollary 2.4]{BGP}, the Jordan--H\"{o}lder factors of the left $F[\mathit{GL}_2(F)]$-module $\Delta_n^\star$
	are
	\[
	\left\{L^\star_{m,\frac{m-n}{2}}\ \middle|\ m\in n-2\cE(n)\right\}.
	\]
	Here we explain why the determinant twist is different from what is claimed in \cite[Corollary 2.4]{BGP}, using the notation of \cite{BGP}.
	In the proof of that corollary, the application of \cite[Theorem 10.1.8]{Bon} is valid after restricting representations via $\varphi_2$.
	However, the normalization given at the beginning of \cite[\S2.1]{BGP} shows that 
	for any $z\in F^\times$, the element $\begin{pmatrix}
		z& 0\\0&z
	\end{pmatrix}\in \mathit{GL}_2(F)$ acts on $\Delta_{k'}\otimes \det^{m'}$ as $z^{-k'+2m'}$ and on $\Delta_k$
	as $z^{-k}$, so that $m'=\frac{k'-k}{2}$.

	Taking the dual, the Jordan--H\"{o}lder factors of the left $F[\mathit{GL}_2(F)]$-module $(\Delta_n^\star)^\vee$
	are
	\[
	\left\{\left(L^\star_{m,\frac{m-n}{2}}\right)^\vee\ \middle|\ m\in n-2\cE(n)\right\}.
	\]
	By Lemma \ref{LemLstarDual}, this equals
	\[
	\left\{W^\star_{m,\frac{n-m}{2}}\ \middle|\ m\in n-2\cE(n)\right\}.
	\]
\end{proof}



\section{$U_t$-action of level $\Gamma_1(t)$ and the key symmetry}\label{SecUt}

Let $k\geq 2$ be an integer and let $n=k-2$. Then Bandini and Valentino give an explicit formula for the action of $U_t$ on $S_k(\Gamma_1(t))$ with respect to a 
basis $\mathbf{c}_0,\ldots,\mathbf{c}_n$ \cite[(4.1)]{BV0}. 
With our normalization, their formula reads as follows (see also \cite[(2.1)]{Ha_GVt}):
\begin{equation}\label{EqnBV}
	\begin{split}
		U_t(\mathbf{c}_j)=&(-t)^j\tbinom{n-j}{j}\mathbf{c}_j\\
		&-t^j\sum_{h\in \bZ, h\neq 0}\left\{\tbinom{n-j-h(q-1)}{-h(q-1)}+(-1)^{j+1}\tbinom{n-j-h(q-1)}{j} \right\}\mathbf{c}_{j+h(q-1)}.
	\end{split}
\end{equation}
Here we use the convention $\binom{c}{d}=0$ unless $0\leq d\leq c$. 
Terms involving $\mathbf{c}_{j+h(q-1)}$ are omitted if $j+h(q-1)\notin 
\{0,\ldots,n\}$. 

Let $U_n\in M_{n+1}(A)$ be the representing matrix of $U_t$ with respect to this basis. We can write
\[
U_n=B_n D_n,\quad B_n\in M_{n+1}(\bF_p),\quad D_n=\diag(1,t,\ldots,t^n).
\]

\begin{lem}\label{LemEntriesOfBn}
	For any $i,j\in \{0,\ldots,n\}$, we have
	\[
	(B_n)_{i,j}=\left\{
	\begin{array}{ll}
		\delta_{i,j}-\binom{n-i}{j-i}+(-1)^j\binom{n-i}{j}& (i\equiv j\bmod q-1),\\
		0 & (i\not\equiv j\bmod q-1).
		\end{array}
	\right.
	\]
\end{lem}
\begin{proof}
	Suppose $i\equiv j\bmod q-1$. If $i\neq j$, then on the right-hand side of (\ref{EqnBV}), the only contribution to $(B_n)_{i,j}$ is the term for $h=\frac{i-j}
	{q-1}
	$. Then $i=j+h(q-1)$ and 
	\[
	(B_n)_{i,j}=-\left(\binom{n-i}{j-i}+(-1)^{j+1}\binom{n-i}{j}\right)=-\binom{n-i}{j-i}+(-1)^{j}\binom{n-i}{j}.
	\]
	If $i=j$, then (\ref{EqnBV}) shows $(B_n)_{i,i}=(-1)^i\binom{n-i}{i}$, whereas 
	\[
	1-\binom{n-i}{0}+(-1)^i\binom{n-i}{i}=(-1)^i\binom{n-i}{i}.
	\]
	This settles the case of $i\equiv j\bmod q-1$. If $i\not\equiv j\bmod q-1$, then (\ref{EqnBV}) yields $(B_n)_{i,j}=0$.
	This concludes the proof.
\end{proof}

Let 
\[
\tau=\begin{pmatrix}
	1&0\\0&t
\end{pmatrix}.
\]
By (\ref{EqnActionDualStar}), we have
\begin{equation}\label{EqnDnGrpRing}
\rho_n^\star(\tau)=D_n.
\end{equation}
For any $b\in \bF_q$ and $c\in \bF_q^\times$, let
\[
u_b=\begin{pmatrix}
	1&b\\0&1
\end{pmatrix},\quad v_c=\begin{pmatrix}
1&-c\\c^{-1}&0
\end{pmatrix}.
\]
Define elements $\cU,\cV,\cB\in F[\mathit{GL}_2(F)]$ by
\[
\cU=\sum_{b\in \bF_q}[u_b],\quad \cV=-\sum_{c\in \bF_q^\times}[v_c],\quad \cB=\cU+\cV.
\]
We have
\begin{equation}\label{EqnTSquareIs0}
	\cU^2=\sum_{b,b'\in \bF_q}[u_b][u_{b'}]=\sum_{b,b'\in \bF_q}[u_{b+b'}]=q\sum_{b\in \bF_q}[u_b]=0.
\end{equation}

\begin{lem}\label{LemBUGrpRing}
	\[
	B_n=\rho_n^\star(\cB),\quad U_n=\rho_n^\star(\cB \tau).
	\]
\end{lem}
\begin{proof}
	By $U_n=B_nD_n$ and (\ref{EqnDnGrpRing}), it suffices to show the first equality.
	
	Take any $i,j\in \{0,\ldots,n\}$. By (\ref{EqnActionDualStar}), the $(i,j)$-entry $\rho^\star_n(\cU)_{i,j}$ of $\rho_n^\star(\cU)$ is equal to the coefficient 
	of $X^{n-j}Y^{j}$ in
	\[
	\sum_{b\in \bF_q} (X+bY)^{n-i}Y^{i}=X^{n-i}Y^i+\sum_{b\in \bF_q^\times} (X+bY)^{n-i}Y^{i}.
	\]
	It is zero unless $q-1\mid j-i$, in which case it equals
	\[
	\delta_{i,j}+\sum_{b\in \bF_q^\times }b^{j-i} \binom{n-i}{j-i}=\delta_{i,j}+(q-1)\binom{n-i}{j-i}=\delta_{i,j}-\binom{n-i}{j-i}.
	\]
	
	On the other hand, $\rho_n^\star(\cV)_{i,j}$ is equal to the coefficient of $X^{n-j}Y^{j}$ in
	\[
	-\sum_{c\in \bF_q^\times} (X-cY)^{n-i}(c^{-1}X)^{i}.
	\]
	The coefficient equals
	\[
	-\sum_{c\in \bF_q^\times}\binom{n-i}{j}(-c)^jc^{-i}=-(-1)^j\binom{n-i}{j} \sum_{c\in \bF_q^\times}c^{j-i}.
	\]
	It is zero unless $q-1\mid j-i$, in which case it equals 
	\[
	(-1)^j\binom{n-i}{j}.
	\]
	Now the result follows from Lemma \ref{LemEntriesOfBn}.
\end{proof}

\begin{lem}\label{LemQTTransposeDifference}
	\[
	 {}^t\!\cV=\cV,\quad \cU-\cV\cU={}^t\!\cU\cU.
	\]
\end{lem}

\begin{proof}
	The first equality holds since ${}^t\! v_c=v_{-c^{-1}}$ and the map
	\[
	\bF_q^\times\to \bF_q^\times,\quad c\mapsto -c^{-1}
	\]
	is a bijection.
	
	Let us show the second equality.
	For any $b\in \bF_q$, let
	\[
	l_b=\begin{pmatrix}
		1&0\\b&1
	\end{pmatrix}={}^t\! u_b.
	\]
	Then $v_c=l_{c^{-1}}u_{-c}$ for any $c\in \bF_q^\times$.
	Since $[u_b] \cU=\cU$ for any $b\in \bF_q$, we have
	\[
	\begin{aligned}
		\cV\cU&=-\sum_{c\in \bF_q^\times} [v_c] \cU=-\sum_{c\in \bF_q^\times} [l_{c^{-1}}][u_{-c}]\cU\\
		&=-\sum_{c\in \bF_q^\times} [l_{c^{-1}}]\cU=-\sum_{c\in \bF_q^\times} [l_{c}]\cU=-({}^t\!\cU-1)\cU,
	\end{aligned}
	\]
	where $1=[\id]$ is the unit element of $F[\mathit{GL}_2(F)]$.
	This yields $\cU-\cV\cU={}^t\!\cU\cU$.
\end{proof}

Define $\cP\in F[\mathit{GL}_2(F)]$ by
\[
\cP=1-\cU.
\]

\begin{prop}\label{PropKeySymmetry}
	\[
	{}^t\!(\cB\cP)=\cB\cP.
	\]
\end{prop}
\begin{proof}
	By (\ref{EqnTSquareIs0}), we have
	\[
	\cB\cP=(\cU+\cV)(1-\cU)=(\cU-\cV\cU)+\cV.
	\]
	By Lemma \ref{LemQTTransposeDifference}, we obtain $\cB\cP=\cV+{}^t\!\cU\cU$ and ${}^t\!(\cB\cP)=\cB\cP$.
\end{proof}

\begin{lem}\label{LemBPCommutesWithY}
	For any $d\in \bF_q^\times$, we have the equalities 
	\[
	\cB y_d=y_d\cB,\quad \cP y_d=y_d\cP,\quad \tau y_d=y_d\tau
	\]
	in $F[\mathit{GL}_2(F)]$.
\end{lem}
\begin{proof}
	We have
	\[
	y_d^{-1}u_b y_d=u_{bd},\quad y_d^{-1}v_c y_d=v_{cd}.
	\]
	Since multiplication by $d$ induces bijections of both $\bF_q$ and $\bF_q^\times$, 
	it follows that $y_d$ commutes with $\cU$ and $\cV$.
	Hence we obtain the first two equalities. The last one is clear.
\end{proof}

For any $m\in \bZ_{\geq 0}$, $a\in \bZ$ and $e\in \{0,\ldots,q-2\}$, let
\[
\cB_{m,a,e}:=\rho_{m,a,e}^\star(\cB),\quad \cP_{m,a,e}:=\rho_{m,a,e}^\star(\cP),
\]
which are well-defined by Lemma \ref{LemBPCommutesWithY}.

\begin{lem}\label{LemBPentries}
	Every entry of $\cB_{m,a,e}$ and $\cP_{m,a,e}$ is an element of $\bF_q$. Moreover, $\cP_{m,a,e}$ is upper unitriangular. 
\end{lem}
\begin{proof}
	Recall that 
	\[
	\cB=\cU+\cV=\sum_{b\in \bF_q}[u_b]-\sum_{c\in \bF_q^\times}[v_c],\quad \cP=1-\cU=1-\sum_{b\in \bF_q}[u_b].
	\]
	Since $u_b,v_c\in \mathit{GL}_2(\bF_q)$, the definition of the representation $\rho^\star_{m,a}$ and (\ref{EqnActionDualStar}) imply 
	the first assertion. 
	
	For the second assertion,
	write $m=m_0+m_1p+\cdots +m_\mu p^\mu$ as before. For any $\nu\in \{0,\ldots,\mu\}$ and $i,j\in \{0,\ldots, m_\nu\}$, (\ref{EqnActionDualStar}) yields
	\[
	(u_b\star (X^{m_\nu-j}Y^j)^\vee)(X^{m_\nu-i}Y^i)=(X^{m_\nu-j}Y^j)^\vee((X+bY)^{m_\nu-i}Y^i).
	\]
	This is zero unless $j\geq i$, in which case it equals $\binom{m_\nu-i}{j-i}b^{j-i}$. 
	Note that, for any $\underline{i}=(i_\nu)_\nu$ and $\underline{j}=(j_\nu)_\nu$, if $j_\nu\geq i_\nu$ for all $\nu$, then the definition of the
	colexicographic order yields $\underline{j}\geq \underline{i}$.
	Since $\det(u_b)=1$, this implies that $\rho^\star_{m,a}(u_b)$ is upper unitriangular. 
	Since any diagonal entry of $\rho^\star_{m,a}(\cP)$ equals $1-\sum_{b\in \bF_q}1=1$, it is also upper unitriangular.
	Then its principal submatrix $\cP_{m,a,e}=\rho^\star_{m,a}(\cP)_{J_e}$ for $J_e$ is 
	also
	upper unitriangular. 
\end{proof}

\begin{prop}\label{PropSymmetricW}
	\[
	{}^t\!(H_{m,J_e}\cB_{m,a,e}\cP_{m,a,e})=H_{m,J_e}\cB_{m,a,e}\cP_{m,a,e}.
	\]
\end{prop}
\begin{proof}
	Since $H_{m,J_e}$ is a diagonal matrix, we have
	\[
	{}^t\!(H_{m,J_e}\cB_{m,a,e}\cP_{m,a,e})={}^t\!(H_{m,J_e}\rho_{m,a,e}^\star(\cB\cP))={}^t\!(\rho_{m,a,e}^\star(\cB\cP))H_{m,J_e}.
	\]
	Then (\ref{EqnSeldual-e-part}) yields
	\[
	{}^t\!(\rho_{m,a,e}^\star(\cB\cP))H_{m,J_e}=H_{m,J_e} \rho_{m,a,e}^\star({}^t\!(\cB\cP)).
	\]
	By Proposition \ref{PropKeySymmetry}, we obtain
	\[
	H_{m,J_e} \rho_{m,a,e}^\star({}^t\!(\cB\cP))=H_{m,J_e} \rho_{m,a,e}^\star(\cB\cP)=H_{m,J_e}\cB_{m,a,e}\cP_{m,a,e}.
	\]
\end{proof}

\begin{cor}\label{CorBSymRankProfile}
	The rank profile matrix of $\cB_{m,a,e}$ is symmetric.
\end{cor}
\begin{proof}
	Since $H_{m,J_e}$ is an invertible diagonal matrix with entries in $\bF_q$, Lemma \ref{LemBPentries} and Proposition \ref{PropSymmetricW} show that 
	the matrix $\cB_{m,a,e}$ satisfies the assumptions of Lemma \ref{LemSymmRankProfile} for $K=\bF_q$.
	Thus, that lemma yields the corollary.
\end{proof}



\section{Proof of the main theorem}\label{SecPfThmMain}

\begin{thm}\label{ThmMain}
	Let $\alpha$ be any finite slope of $U_t$ acting on $S_k(\Gamma_1(t^r))$.
	\begin{enumerate}
		\item\label{ThmMain_HI} $\alpha\in \frac{1}{2}\bZ$. If $q$ is odd, then $\alpha\in \bZ$.
		\item\label{ThmMain_Bound} $\alpha\leq k-2$.
	\end{enumerate}
\end{thm}
\begin{proof}
	By \cite[Theorem 1.1 (1)]{Ha_DMDia}, applied with $\frn=1$ and $\wp=t$, the sets of finite slopes of $U_t$ on $S_k(\Gamma_1(t^r))$ 
	and $S_k(\Gamma_1(t))$ coincide. Therefore, it suffices to treat the case $r=1$.
	
	Let $n=k-2$.
	By Lemma \ref{LemBUGrpRing}, it suffices to show that the finite slopes of the matrix $\rho_n^\star(\cB \tau)$ satisfy the two asserted properties. Note that 
	the multiset of slopes of this matrix is the sum of multisets of the slopes of representing matrices of $\cB \tau\in F[\mathit{GL}_2(F)]
	$ 
	acting 
	on Jordan--H\"{o}lder 
	factors of $
	(\Delta_n^\star)^\vee$. Thus, by Proposition \ref{PropJH},
	it is enough to show these properties for the finite slopes of $\rho_{m,\frac{n-m}{2}}^\star(\cB \tau)$ for any $m\in n-2\cE(n)$.
	
	Fix $m\in n-2\cE(n)$. Let $m=m_0+m_1 p+\cdots +m_\mu p^\mu$ be the base-$p$ expansion of $m$. Let $a=\frac{n-m}{2}$.
	By (\ref{Eqnn-mhalfWellDef}), we have $a\in \bZ_{\geq 0}$.
	For any $e\in \{0,\ldots, q-2\}$, consider the $F$-subspace $W^\star_{m,a,e}$ of $W^\star_{m,a}$ as before.
	By Lemma \ref{LemBPCommutesWithY}, the action of $\cB \tau$ stabilizes each $W^\star_{m,a,e}$.  
	Note that we have
	\[
	\det(u_b)=\det(v_c)=1,\quad \det(\tau)=t.
	\]
	Since the $F$-vector space $W^\star_{m,a}$ is the direct sum of the $W^\star_{m,a,e}$, it is enough to show the two asserted properties
	for the finite slopes of 
	\begin{equation}\label{EqnSlopesAzero}
	\rho_{m,a,e}^\star(\cB \tau)=t^{a}\rho_{m,0,e}^\star(\cB \tau).
	\end{equation}
	If $J_e=\emptyset$, then this block is empty and contributes no eigenvalues. Hence, we may assume $J_e\neq\emptyset$.

	From (\ref{EqnActionDualStar}), we see that the action of $\tau$ on $W^\star_{m,0,e}$ is given by
	\[
	f^\wedge_{\underline{i}}\mapsto t^{w(\underline{i})}f_{\underline{i}}^\wedge,\quad \underline{i}=(i_\nu)_{\nu\in \{0,\ldots,\mu\}}\in J_e,
	\]
	where $w(\underline{i})=\sum_{\nu=0}^\mu i_\nu p^\nu\in \bZ_{\geq 0}$ as before.
	Let $d=\dim_F(W^\star_{m,0,e})\geq 1$.
	We can write
	\[
	D_{m,0,e}:=\rho^\star_{m,0,e}(\tau)=\diag(t^{a_0},\ldots,t^{a_{d-1}}),
	\]
	where for any $j\in\{0,\ldots,d-1\}$, we have $a_j=w(\underline{i})$ for the $(j+1)$-st smallest element $\underline{i}\in J_e$
	in the colexicographic order. Then (\ref{EqnCoLexVSWeight}) yields
	\[
	a_0<\cdots <a_{d-1}.
	\]
	
	By Lemma \ref{LemBPentries}, we have $\cB_{m,0,e}\in M_d(\bF_q)$. 
	By Corollary \ref{CorBSymRankProfile}, the matrix $\cB_{m,0,e}$ has a symmetric rank profile matrix.
	Then the matrices $D_{m,0,e}$ and $\cB_{m,0,e}$ 
	satisfy the assumptions of Theorem \ref{ThmSymSlopes} for $K=\bF_q$. Therefore,
	by Theorem \ref{ThmSymSlopes} and (\ref{EqnSlopesAzero}), every finite slope of $\rho_{m,a,e}^\star(\cB \tau)$
	is equal to 
	\begin{equation}\label{EqnSlopeValue}
	a+\frac{1}{2}(w(\underline{i})+w(\underline{j}))
	\end{equation}
	for some $\underline{i},\underline{j}\in J_e$.
	
	By the definition of the set $J_e$, we can write
	\[
	w(\underline{i})=N(q-1)+e,\quad w(\underline{j})=N'(q-1)+e
	\]
	for some $N,N'\in \bZ$.
	Thus (\ref{EqnSlopeValue}) equals
	\[
	a+e+(N+N')\left(\frac{q-1}{2}\right),
	\]
	which belongs to $\frac{1}{2}\bZ$, and is an integer if $q$ is odd. This proves (\ref{ThmMain_HI}).
	Moreover, (\ref{EqnRangeOfWeight}) and (\ref{Eqnn-mhalfWellDef}) yield
	\[
	\begin{aligned}
		a+\frac{1}{2}(w(\underline{i})+w(\underline{j}))\leq \frac{n-m}{2}+m=\frac{n+m}{2}\leq n=k-2.
	\end{aligned}
	\]
	This proves (\ref{ThmMain_Bound}).
\end{proof}



\end{document}